\documentclass[a4paper, 12pt, leqno]{amsart}
\usepackage{mathrsfs}
\usepackage{amsmath}
\usepackage{amscd}
\usepackage{amssymb}
\usepackage{amsthm}

\numberwithin{equation}{section}

\def\BN{\mathbb{N}}

\def\BZ{{\mathbb{Z}}}

\def\cala{\mathcal{A}}

\begin{document}
\title[Lusztig's $a$-function for Coxeter groups]
{Lusztig's $a$-function for Coxeter groups of rank 3}

\author{ Peipei Zhou}

\begin{abstract}
We show that  Lusztig's $a$-function of a Coxeter group is bounded
if the rank of the Coxeter group is 3.
\end{abstract}
\maketitle \setcounter{section}{-1}
\section{Introduction}
\def\tgt{\tilde G_2}
\label{sec:Intro} In [L2] Lusztig defined
$a$-function for a Coxeter group and showed that $a$-function is bounded
for affine Weyl groups. This boundness plays an important role in studying
cells of affine Weyl groups. In [X], Xi showed that the
$a$-function is bounded for Coxeter groups with complete Coxeter graph.
He also gave some interesting applications of the boundness on cells of the Coxeter groups.
In this paper, we show that Lusztig's $a$-function of a Coxeter group is bounded if the rank of the Coxeter group is 3.
The present work was motivated by a question posed by Prof. Xi in his paper [X]. The author would like to thank Prof. Xi for his help in
dealing with the problems in writing the paper.

\medskip

\section{Preliminaries}

\noindent{\bf 1.1.} We first recall some known facts,
and refer to [KL, L2, L3,X] for more details.
 Let $(W,S)$ be a Coxeter group.  Denote $l$
the length function and $\le$  the Bruhat order of $W$. The
neutral element of $W$ will be denoted by $e$.

Let $q$ be an indeterminate. The Hecke algebra $H$ of $(W,S)$ is a
free $\cala=\BZ[q^{\frac12},q^{-\frac12}]$-module with a basis
$T_w,\ w\in W$ and the multiplication relations are
$(T_s-q)(T_s+1)=0$ if $s$ is in $S$, $T_wT_u=T_{wu}$ if
$l(wu)=l(w)+l(u)$.

\def\tt{\tilde T}

For any $w\in W$ set $\tilde T_w=q^{-\frac{l(w)}2}T_w$. For any
$w,u\in W$, write $$\tilde T_w\tilde T_u=\sum_{v\in
W}f_{w,u,v}\tilde T_v,\qquad f_{w,u,v}\in \cala.$$ The following fact is
known and implicit in [L2, 8.3], see also [X] 1.1.(a).

\medskip

 \noindent{(a)} For any $w,u,v\in W$, $f_{w,u,v}\in\cala $ is a polynomial in
$q^{\frac12}-q^{-\frac12}$ with non-negative coefficients  and
$f_{w,u,v}=f_{u,v^{-1},w^{-1}}=f_{v^{-1},w,u^{-1}}$. Its degree is
less than or equal to min$\{l(w),l(u), l(v)\}$.

\medskip

For any $w,u,v$ in $W$, we shall regard $f_{w,u,v}$ as a polynomial
in $\xi=q^{\frac12}-q^{-\frac12}$. The following  fact is due to
Lusztig [L3, 1.1 (c)].

\medskip

\noindent(b) For any $w,u,v$ in $W$ we have
$f_{w,u,v}=f_{u^{-1},w^{-1},v^{-1}}.$

\medskip
We shall need the following facts.

\noindent(c) Let $(W,S)$ be a Coxeter group and $I$ is
a subset of $S$. The following conditions are equivalent.

(1) The subgroup $W_I$ of $W$ generated by $I$ is finite.

(2) There exists an element $w$ of $W$ such that $sw\le w$ for all
$s$ in $I$.

(3) There exists an element $w$ of $W$ such that $w\le ws$ for all
$s$ in $I$.

\medskip

As usual, we set $L(w)=\{s\in S\,|\, sw\le w\}$ and $R(w)=\{s\in S\,|\, ws\le
w\}$ for any $w\in W$.

\medskip

\noindent(d) Let $w$ be in $W$ and $I$ is a subset of
$L(w)$ (resp. $R(w)$). Then $l(w_Iw)+l(w_I)=l(w)$ (resp.
$l(ww_I)+l(w_I)=l(w)$), here $w_I$ is the longest element of $W_I$.

\medskip

\noindent{\bf 1.2.} For any $y,w\in W$, let $P_{y,w}$ be the
Kazhdan-Lusztig polynomial. Then all the elements
$C_w=q^{-\frac{l(w)}2}\sum_{y\le w}P_{y,w}T_y$, $w\in W$, form a
Kazhdan-Lusztig basis of $H$. It is known that
$P_{y,w}=\mu(y,w)q^{\frac12(l(w)-l(y)-1)}$ +lower degree terms if
$y<w$ and $P_{w,w}=1$.

For any $w,u$ in $W$, Write $$C_wC_u=\sum_{v\in W}h_{w,u,v}C_v,\
h_{w,u,v}\in\cala.$$ Following [L2], for any $v\in W$ we define
$$a(v)=\max\{i\in\BN\,|\, i=\text{deg}h_{w,u,v}, \ w,u\in W\},$$
here the degree is in terms of $q^{\frac12}$. Since $h_{w,u,v}$ is a
polynomial in $q^{\frac12}+q^{-\frac12}$, we have $a(v)\ge 0$.

We are interested in the bound of the function $a:W\to\BN$. Clearly,
$a$ is bounded if $W$ is finite. The following fact is known (see
[L3]).

\medskip
The $a$-function is bounded by a constant $c$ if and
only if deg$f_{w,u,v}\le c$ for any $w,u,v\in W$.

Lusztig showed that for an affine Weyl group the $a$-function is
bounded by the length of the longest element of the corresponding
Weyl group. This fact is important in studying cells in affine Weyl
groups. One consequence is that an affine Weyl group has a lowest
two-sided cell [S1]. In general, Xi showed that the lowest two-sided
cells exists for a Coxeter group with bounded $a$-function. (see [X,1.5])
\medskip

\section{Coxeter groups of rank 3 }

\def\st{\stackrel}
\def\sc{\scriptstyle}

In this section $(W,S)$ is a infinite Coxeter group of rank 3. Let $S=\{r,s,t\}$, we shall assume that $tr=rt$.
By 1.1.(c), for $w\in W$, both $R(w)$ and $L(w)$ contain at most 2 elements.
Let $|R(w)|$ (resp. $|L(w)|$) denote the number of elements in the set $R(w)$ (resp. $L(w)$). Let $m_{sr}$ (resp. $m_{st}$) denote the order of $sr$ (resp. $st$).
Let $w_{sr}$ (resp. $w_{st}$) denote the longest element in the parabolic subgroup generated by $s, r$ (resp. $s, t$).

\medskip

 \noindent{\bf Theorem 2.1.}
 Let $(W,S)$ be a Coxeter group of rank 3 and assume that $rt=tr$, $S=\{r,s,t\}$. Then Lusztig's $a$-function on $W$ is bounded by the
length of the longest element of certain finite parabolic subgroups
of $W$ in the following two cases:

(a) $m_{sr}\geq 7$ and $m_{st}=3 $.

(b) $m_{sr}\geq 5$ and $m_{st}\geq 4 $.

The remaining of this paper is devoted to a proof of the theorem.

\medskip
 In section 3, we will deal with the case (a).

 In section 4, we will deal with the case (b).

\medskip
\textbf{The notation}, if $w=(w_{1})(w_{2})\cdots (w_i)$, means $w=w_{1}w_{2}\cdots w_i$ and $l(w)=l(w_{1})+l(w_{2})+\cdots +l(w_i)$.
\medskip

The strong exchange condition will be need frequently in the proof, so we recall it.
\medskip

 \noindent{\bf Strong exchange condition.}
Let $(W, S)$ be a Coxeter group. Let $w=s_1\cdots s_r (s_i\in S)$, not necessarily a reduced expression. Suppose $t\in \bigcup_{w\in W}wSw^{-1}$, satisfies $l(wt)<l(w)$. Then there is an index $i$ for which $wt=s_1\cdots \hat{s_i}\cdots s_r$ (omitting $s_i$). If the expression for $w$ is
 reduced, then $i$ is unique.

\section{The case $m_{sr}\geq 7$ and $m_{st}=3 $}
Since $m_{st}=3$, $m_{sr}\geq 7$, $w_{st}=sts=tst$ and $l(w_{sr})\geq 7$.

\medskip

\noindent{\bf Lemma 3.1.}
There is no element $w$ in $W$ such that $w=(w_1)(st)=(w_2)(sr)$ .

Proof. We use induction on $l(w)$. When $l(w)=0,1,2,3$, the lemma is
clear. Now assume that the lemma is true for $u$ with
$l(u)\leq l(w)-1$. Since $r,t\in R(w)$, by 1.1.(d),
$w=(w_3)(rt)$ for some $w_3\in W$.
So we get $w_1s=w_3r,w_2s=w_3t$. By 1.1.(d),
$w_1s=w_3r=(w_4)(w_{sr})$ for some $w_4\in W$,
$w_2s=w_3t=(w_5)(w_{st})$ for some $w_5\in W$. Since $m_{sr}\geq 7$,
we have $\widetilde{w_4}\in W$, such that $w_3=(w_4)(w_{sr}r)=$ $(\widetilde{w_4})(srsrs)=$
$(w_5)(w_{st}t)=(w_5)(ts)$. Then
there exists $w_6,w_7 \in W$, such that $w_7=(\widetilde{w_4})(sr)=(w_6)(st)$.
By induction hypothesis, $w_7 $ does not exist, hence $w$ does not exist. The
lemma is proved.

\medskip
\noindent{\bf Corollary 3.2.}
There is no element $w$ in $W$ such that $w=(w_1)(srs)=(w_2)(t)$ .

Proof.Assume that $w$ exists, by 1.1.(d), there exists $w_3\in W$ , such that $w=(w_3)(w_{st})$, hence $(w_1)(sr)=(w_3)(st)$, which contradicts Lemma 3.1.

\medskip
\noindent{\bf Corollary 3.3.}
There is no element $w$ in $W$ such that $w=(w_{1})(srsr)=(w_{2})(t)$ .

Proof. Assume that $w$ exists, there exists $w_3\in W$, such that $w=(w_3)(tr)$, hence $(w_1)(srs)=(w_3)(t)$, which contradicts Corollary 3.2.

\medskip

\noindent{\bf Lemma 3.4.}
There is no element $w$ in $W$ such that $w=(w_{1})(ts)=(w_{2})(r)$ .

Proof. Assume that $w$ exists, there exists $w_3\in W$, such that $w=(w_3)(w_{sr})$, hence $(w_1)(t)=(w_3)(w_{sr}s)$, which contradicts Corollary 3.2.

\medskip

\noindent{\bf  Lemma 3.5.}
 Let $x, y$ be elements in $W$, and $w$ be an element in the parabolic subgroup $W_{sr}$ generated by $r, s$.
 Assume that $l(w)\ge 5$ and  $r, s\notin R(x)\cup L(y)$. Then

\noindent (a) $xwy=(x)(w)(y)$, i.e., $l(xwy)=l(x)+l(w)+l(y)$.

\noindent (b) $R(xwy)=R(wy)$.

\noindent (c) $L(xwy)=L(xw)$.

Proof: It is clear that $xw=(x)(w)$, and $wy=(w)(y)$.
 Note that (b) and (c) are equivalent. We use induction on $l(y)$ to prove (a) and (b).
 The case $l(y)=0$ is clear. If $l(y)=1$, then $y=t$. By Corollaries 3.2 and 3.3, we see $xwt=(x)(w)(t)$.
 If $R(wt)$ contains two elements, we must have $R(xwt)=R(wt)$. Now assume that $R(wt)=\{t\}$.
 If $R(xwt)\neq R(wt)$, $R(xwt)=\{r, t\}$, or $\{s,t\}$.
 If $R(xwt)=\{r,t\}$, we have $r\in R(xw)=R(w)$, which contradicts that $R(wt)=\{t\}$.
 If $R(xwt)=\{s,t\}$, we have $xwt=(u)(tst)$, for some $u\in W$. Then $xw=(u)(ts)$,
 so $w=(w_1)(srsrs)$ for some $w_1\in W_{sr}$. By Corollary 3.3, this is impossible.
 Hence $R(xwt)=R(wt)$.

 If $l(y)=2$, then $y=ts$. By what we have proved that $s\notin R(xwt)$, we see that $xwts=(x)(w)(ts)$.
 If $R(wts)$ contains two elements, we must have $R(xwts)=R(wts)$. Now assume that $R(wts)=\{s\}$.
 If $R(xwts)\neq R(wts)$, $R(xwts)=\{s, t\}$, or $\{s,r\}$.
 If $R(xwts)=\{s,t\}$, we have $s\in R(xw)=R(w)$, which contradicts that $R(wts)=\{s\}$.
 If $R(xwts)=\{s,r\}$, we have $xwts=(u)(w_{sr})$, for some $u\in W$. Then $xwt=(u)(w_{sr}s)$,
 since $w_{sr}s=(w_1)(srsrsr)$ for some $w_1\in W_{sr}$, by Corollary 3.3, this is impossible.
 Hence $R(xwts)=R(wts)$.

Now assume that $k\geq 3$. Let $y=y_1y_2\cdots y_k$ be a reduced expression of $y$.
The induction hypothesis says that $R(xwy_1\cdots y_i)=R(wy_1\cdots y_i)$
and $l(xwy_1\cdots y_i)=l(x)+l(w)+i$ for $i\leq k-1$.
We must have $y_{k}\notin R(xwy_1y_2\cdots y_{k-1})$, since $wy=(w)(y)$,
so $xwy=(x)(w)(y)$.

Assume that $|R(xwy_1y_2\cdots y_{k-1})|=2$.
If $R(xwy)$ contains one element, it must be $y_k$,
so $R(xwy)=R(wy)=\{y_k\}$.
When $R(xwy)$ contains two elements,
if $R(xwy_1y_2\cdots y_{k-1})=\{r,s\}$ or $\{t,s\}$, then $y_k=t$ or $r$,
and $R(xwy)=R(wy)=\{t,r\}$.
If $R(xwy_1y_2\cdots y_{k-1})=\{r,t\}$, then $y_k=s$.
When $R(wy)$ contains two elements, we must have $R(xwy)=R(wy)$.
When $R(wy)=\{s\}$, we need to show that $R(xwy)=\{s\}$.
Otherwise $R(xwy)=\{s,r\}$, or $\{s,t\}$.
By Lemma 3.4, $r\notin R(xwy)$, then $R(xwy)=\{s,t\}$.
By 1.1.(d), we have $xwy_1\cdots y_k=(u_1)(sts)$, for some $u_1\in W$.
Then $xwy_1\cdots y_{k-2}y_{k-1}=(u_1)(st)$.

We discuss it in the following three conditions:

(1) $\{y_{k-2}, y_{k-1}\}=\{t, r\}$,
under this condition, $xwy_1\cdots y_{k-3}r=(u_1)(s)$. Hence $r, s\in R(wy_1\cdots y_3r)$.
By 1.1.(d), there exists $u_2\in W$ and $wy_1\cdots y_{k-3}r=(u_2)(w_{sr})$.
Hence $R(wy_1\cdots y_{k-3}rts)=\{s, t\}$, which contradicts to $R(wy)=\{s\}$.

(2)$\{y_{k-2}, y_{k-1}\}=\{s, r\}$.
This is impossible since under this condition, $y_{k-1}=r, y_{k-2}=s$,
from the above, we have $(xwy_1\cdots y_{k_3})(sr)=(u_1)(st)$, which contradicts Lemma 3.1.

(3)$\{y_{k-2}, y_{k-1}\}=\{s, t\}$.
Under this condition, it is easy to see that
$R(wy)=\{s, t\}$, which contradicts to $R(wy)=\{s\}$.

Next assume that $|R(xwy_1y_2\cdots y_{k-1})|=1$,
so $R(xwy_1y_2\cdots y_{k-1})$ $=\{y_{k-1}\}$.
If $R(wy)$ contains two elements, we must have $R(xwy)=R(wy)$.
If $R(wy)$ contains one elements, we must have $R(wy)=\{y_k\}$.
We need to prove that $R(xwy)=\{y_k\}$.
Assume that $R(xwy)\supsetneqq R(wy)$.
 When $R(xwy)=\{t, r\}$,
it is easy to see that $\{y_{k-1}, y_k\}=\{t,r\}$, so $R(wy)=\{t,r\}$, which contradicts to $R(wy)=\{y_k\}$.

When $R(xwy)=\{t, s\}$, it is easy to see that $y_{k-1}=t$ and $y_k=s$
(or $y_{k-1}=s$ and $y_k=t$), then $t\in R(wy_1\cdots y_{k-2})$ (or $s\in R(wy_1\cdots y_{k-2})$).
So $R(wy)=\{s,t\}$, which is a contradiction.

When $R(xwy)=\{s,r\} $, then $\{y_{k-1},y_k\}=\{s,r\}$.
By 1.1.(d), there exists $u_1\in W$, such that $xwy=(u_1)(w_{sr})$,
write $wy_1\cdots y_k=$
$wy_1\cdots y_i$ $s^a(rs)^br^c$, here $i$ is minimal, such that $R(y_1\cdots y_i)=\{y_i\}=\{t\}$,
$a, c=0$ or 1, $b\geq 0$.

Obviously $a+2b+c< m_{sr}$. Write $u_{sr}=w_{sr}r^c(rs)^{-b}s^a$, then $l(u_{sr})>0$, $xwy_1\cdots y_i=(u_1)(u_{sr})$.
Assume that $R(wy_1\cdots y_i)=R(y_1\cdots y_i)=\{y_i\}$, then $i=0$. So $y\in W_{sr}$, $W_{sr}$ is the parabolic subgroup generated by $s,r$. Then it contradicts to
$s,r\notin L(y)$.

Next assume that $R(wy_1\cdots y_i)\supsetneqq \{y_i\}$.

Only consider $a+2b+c\leq m_{sr}-2$, since when $a+2b+c=m_{sr}-1$, $R(wy)=\{s,r\}$, which contradicts the assumption.

By Corollary 3.2, if $R(u_{sr})=\{s\}$, then $l(u_{sr})\leq 2 $, we must have $u_{sr}=rs$. Hence $R(wy_1\cdots y_i)=\{t,s\}$. If $i$ is large enough, suppose $i\geq 6$, we will show this is impossible.

By the assumption and easy calculation, we get $y_i=t$, $y_{i-1}=s$, $y_{i-2}=r$, $y_{i-3}=s$.
Next we shall deal with the following two conditions:

1) $y_{i-4}=t$, then $y_{i-5}=r$. Hence $R(wy_1\cdots y_{i-5})=\{s,r\}$. By 1.1.(d), there exists $u_2\in W$, such that $xwy_1\cdots y_{i-5}tsrst$ $=(xu_2)$ $(w_{sr}tsrst)$ $=(u_1)(rs)$. By what we have proved already and easy calculation, we see that
there exists $u_3\in W$, such that $(xu_2)(w_{sr}s)(t)=(u_3)(srsr)$, which contradicts Corollary 3.3.

2) $y_{i-4}=r$, and $y_{i-5}=s$, since $y_{i-5}=t$ is as same as condition 1).
Hence $t\in R((xwy_1\cdots y_{i-6})(srsr))$.
But by Corollary 3.3, this is a contradiction .

When $i\leq 5$, there are two cases which satisfy the assumption.

1) $i=2$, $y_1=t$, $y_2=s$, however it contradicts to fact that $y_i=t$, since $i=2$.

2) $i=5$, $y_ 5=t$, $y_4=s$, $y_3=r$, $y_2=s$, $y_1=t$, and $s\in R(w)$. We will show that $xwtsrst=(u_1)(rs)$ is impossible.
Otherwise, $(xws)(tsrst)=(u_1)(r)$, by 1.1.(d), there is $u_3\in W$, such that $(xws)(tsrs)=(u_3)(r)$, by 1.1.(d) there exists $u_4\in W$, such that $(xws)(t)=(u_4)(rsrs)$, which contradicts Corollary 3.3.

If $R(u_{sr})=\{r\}$, then $l(u_{sr})\leq 3$, $R(wy_1\cdots y_i)=\{r,t\}$.
Hence $r,s \in R(wy_1\cdots y_{i-1})$.
By 1.1.(d), there exists $u_5\in W$, such that $wy_1\cdots y_{i-1}=(u_5)(w_{sr})$. Then we get the formula
$(x)(u_5)(w_{sr})(t)=(u_1)(u_{sr})$, here $u_{sr}=sr$ or $u_{sr}=rsr$, however the formula contradicts Lemma 3.1.

Until now, we see that the lemma is proved.

\medskip
\def\tt{\tilde T}
Recall that $\tt_x\tt_y=\sum_{z\in W}f_{x,y,z}\tt_{z}$. Here $f_{x,y,z}$ is a polynomial in $\xi$, where $\xi=(q^{\frac{1}{2}}-q^{-\frac{1}{2}})$.

\textbf{Definition}: deg $\tt_x\tt_y=$  max$_{ z\in W}$ $\{$ deg $f_{x,y,z}\} $.

\medskip

\noindent{\bf Lemma 3.6.} Let $x,y\in W$. Assume that $s,t\notin R(x)\cup L(y)$, then deg $f_{xsts,y,z}\le 1$ for
all $z$ in $W$.

Proof. Write $y=y_1\cdots y_k$, reduced decomposition. Let $w=sts=tst$. There are two cases to consider.

 Case 1: There is no $x'\in W$, such that $x=(x')(w_{sr}s)$, we claim that $R(xwu)=R(wu)$, with $L(u)=\{r\}$, hence the corollary, $l(xstsy)=l(x)+3+l(y)$.
Hence deg $f_{xsts,y,z}=0$, for all $z\in W$.

We use induction on $l(y)$ to prove the claim. When $l(y)=0, 1, 2$, it is easy to see that
$R(xwy)=R(wy)$. When $k\geq 3$, now  assume that the claim is true for $u\in W$, with $l(u)<k$, $t,s\notin L(u)$.
From the proof of Lemma 3.5, we only need to prove the lemma when $R(wy_1\cdots y_{k-1})=\{y_{k-1}\}$ and $R(xwy)=\{s, r\}$.
It is easy to check that $R(wy)\subset R(xwy)$, when $R(wy)$ contains two elements, we must have $R(wy)= R(xwy)$, nothing needs to prove.
Assume that $R(wy)\subsetneqq R(xwy)$.
It is easy to check that $\{y_{k-1}, y_k\}=\{s, r\}$. Write $xwy=(u_1)(w_{sr})$ by 1.1.(d), $u_1\in W$.
Write $wy=wy_1\cdots y_is^a(rs)^br^c$, $0\leq i\leq k-2$, $i$ minimal such that $R(y_1\cdots y_i)=\{y_i\}=\{t\}$, $a+2b+c<m_{sr}$, then we have $xwy_1\cdots y_i=(u_1)(u_{sr})$, $u_{sr}=w_{sr}r^c(rs)^{-b}$ .
From the proof of Lemma 3.5, we only have to check the case $R(wy_1\cdots y_i)=\{s,t\}$ and $i\leq 5$, $u_{sr}=rs$.
By calculation only $i=3$, $y_1=r$, $y_2=s$, $y_3=t$ satisfies the assumption $R(wy_1y_2y_3)=\{s,t\}$ and $R(y_1y_2y_3)=\{t\}$. However, $xwrst=(u_1)(rs)$ is impossible. Otherwise, we will get $(x)(t)=(u_2)(srsr)$, $u_2\in W$, which contradicts Corollary 3.3.
Hence the claim.

Case 2: When there exists $x'\in W$, such that $x=(x')(w_{sr}s)$, $xsts=(x')(w_{sr})(ts)$.We claim that deg $\tt_{xsts}\tt_y= 1$.
If $l(xstsy)=l(x)+3+l(y)$, then nothing needs to prove.
We calculate $\tilde T_{w_{sr}ts}\tilde T_y$ firstly. Assume that there is an $i<k$, which is minimal, such that $l(w_{sr}tsy_1\cdots y_i)<l(w_{sr}tsy_1\cdots y_{i-1})$.
By strong exchange condition, and $l(stsy)=3+l(y)$, we get $rsy_1\cdots y_{i-1}=sy_1\cdots y_{i-1}y_i$,
$l(rsy_1\cdots y_{i-1})=i+1$.
By 1.1.(d), there exists $u_1\in W$, such that $sy_1\cdots y_{i-1}y_i=(w_{sr})(u_1)$.
Since $l(sy)=l(y)+1$, $rsy_1\cdots y_{i-1}=sy_1\cdots y_{i-1}y_i$, hence $rsy_1\cdots y_{i-1}y_{i+1}\cdots y_k$ is a reduced decomposition, and so is
$sy_1\cdots y_{i-1}y_{i+1}\cdots y_k$.

\begin{eqnarray*}
\tilde T_{w_{sr}ts}\tilde T_y & = &  \tilde T_{w_{sr}ts}\tilde T_{y_1\cdots y_{i-1}}\tilde T_{y_i}\tilde T_{y_{i+1}\cdots y_k} \\
                              & = &  \xi \tilde T_{w_{sr}rt}\tilde T_{sy_1\cdots y_{i-1}y_i}\tilde T_{y_{i+1}\cdots y_k}+
                              \tilde T_{w_{sr}rt}\tilde T_{sy_1\cdots y_{i-1}}\tilde T_{y_{i+1}\cdots y_k}  \\
                              & = &  \xi \tilde T_{w_{sr}rt}\tilde T_{sy}+ \tilde T_{w_{sr}rt}\tilde T_{sy_1\cdots y_{i-1}y_{i+1}\cdots y_k}\\
\end{eqnarray*}
Because $sy_1\cdots y_{i-1}y_i=(w_{sr})(u_1)$, $sy=(w_{sr})(u_1)(y_{i+1}\cdots y_k)$,
$sy_1\cdots y_{i-1}$ $y_{i+1}\cdots y_k=(rw_{sr})(u_1)(y_{i+1}\cdots y_k)$.

Since $l(w_{sr})$, $l(rw_{sr})\geq 5$ and $R(w_{sr}rt)=\{t\}$, $\tilde T_{w_{sr}rt}\tilde T_{sy} =\tilde T_{w_{sr}rtsy}$, by Lemma 3.5.

By Lemma 3.5, $\tilde T_{w_{sr}rt}\tilde T_{sy_1\cdots y_{i-1}y_{i+1}\cdots y_k}=\tilde T_{w_{sr}rtsy_1\cdots y_{i-1}y_{i+1}\cdots y_k}$.

Since $l(w_{sr}r)\geq 6$,
then by Lemma 3.5, $\tilde T_{x'}\tilde T_{w_{sr}rtsy}=\tilde T_{x'w_{sr}rtsy}$,
$\tilde T_{x'}\tilde T_{w_{sr}rtsy_1\cdots y_{i-1}y_{i+1}\cdots y_k}=\tilde T_{x'w_{sr}rtsy_1\cdots y_{i-1}y_{i+1}\cdots y_k}$.

Hence the lemma is proved.

\def\tt{\tilde T}

\medskip

\noindent{\bf  Lemma 3.7.} Let $x,y\in W$. Assume that $t,r\notin R(x)\cup L(y)$, then deg $f_{xtr,y,z}\le 2$ for
all $z$ in $W$.

Proof. There are four cases:

 Case 1: When there is no $x'\in W$, such that $x=(x')(w_{sr}r)$, or $x=(x')(w_{st}t)$, or $x=(x')(w_{sr}sr)$.
Claim that $R(xtru)=R(tru)$, with $t,r\notin L(u)$, $u\in W$. Then we have the corollary, $xtry=(x)(tr)(y)$.
Hence deg $f_{xtr,y,z}=0$. We use induction on $l(y)$ to prove the claim. When $l(y)=0,1,2$, it is easy to see the claim is true.
When $l(y)\geq 3$, write $y=y_1y_2\cdots y_k$, reduced decomposition. Now assume that the claim is true for $u$ with $l(u)<k$, $r,t\notin L(u)$.
By the proof of Lemma 3.5, we only have to prove that when $R(wy_1\cdots y_{k-1})=\{y_{k-1}\}$
and $R(xtry)=\{s, r\}$, $R(try)=\{s, r\}$.
Assume that $R(try)\subsetneqq R(xwy)$. It is easy to see that $\{y_{k-1},y_k\}=\{s,r\}$. Write $xtry=(u_1)(w_{sr})$ by 1.1.(d), for $u_1\in W$.
Write $try=try_1\cdots y_is^a(rs)^br^c$, $0\leq i\leq k-2$, $i$ minimal such that $R(y_1\cdots y_i)=\{y_i\}=\{t\}$, $a+2b+c<m_{sr}$, then we have $xtry_1\cdots y_i=(u_1)(u_{sr})$, where $u_{sr}=w_{sr}r^c(rs)^{-b}$ .
From the proof of Lemma 3.5, we only have to consider the case $R(try_1\cdots y_i)=\{s,t\}$ and $i\leq 5$, $u_{sr}=rs$.

If $i=5$, $y_5=t$, $y_4=s$, $y_3=r$, $y_2=s$,contradicts to  $y_1=s$, and $y_1y_2\cdots y_5$ is a reduced decomposition.

If $i=4$, then $y_4=t$, $y_3=s$, $y_2=r$, $y_1=s$,
this contradicts to the assumption $R(try_1\cdots y_i)=\{s,t\}$.

If $i=3$, then $y_3=t$, $y_2=s$, which contradicts to $sy_2\cdots y_i$ is a reduced decomposition.

It is easy to see that $i=1$ is impossible.

If $i=2$, $y_2=t$, $y_1=s$, which satisfies $R(trst)=\{s,t\}$ and $R(st)=\{t\}$. That is $y_3\cdots y_k=srw_{sr}$, $ty_1\cdots y_k=rstrw_{sr}$. However the equality $xtrst=(u_1)(rs)$ is failed to hold, since there is no $x'\in W$, such that $x=(x')(w_{sr}sr)$.

Hence the claim.

Case 2: When $x=(x')(w_{sr}r)$, then $xtr=(x')(w_{sr})(t)$. If $l(xtry)=l(x)+2+l(y)$, then nothing needs to prove.
First we calculate $\tilde T_{w_{sr}t}\tilde T_y$.  Assume that there is an $i<k$, which is minimal, such that $l(w_{sr}ty_1\cdots y_i)$
$ <l(w_{sr}ty_1\cdots y_{i-1})$.
By strong exchange condition, and $l(try)=2+l(y)$, we get $sty_1\cdots y_{i-1}=ty_1\cdots y_{i-1}y_i$,
$l(sty_1\cdots y_{i-1})=i+1$.

Since $l(try)=l(y)+2$ and $sty_1\cdots y_{i-1}=ty_1\cdots y_{i-1}y_i$, it is easy to see $sty_1\cdots y_{i-1}y_{i+1}\cdots y_k$ and $ty_1\cdots y_{i-1}y_{i+1}\cdots y_k$ are reduced decompositions.
Let $\xi=(q^{\frac{1}{2}}-q^{-\frac{1}{2}})$.
\begin{eqnarray*}
\tilde T_{w_{sr}t}\tilde T_y & = &  \tilde T_{w_{sr}s}\tilde T_{sty_1\cdots y_{i-1}}\tilde T_{y_i}\tilde T_{y_{i+1}\cdots y_k} \\
                              & = &  \xi \tilde T_{w_{sr}s}\tilde T_{ty}+ \tilde T_{w_{sr}s}\tilde T_{ty_1\cdots y_{i-1}y_{i+1}\cdots y_k}
\end{eqnarray*}
We have showed that $ L(ty)=\{s,t\}$.
If there exists $u_1\in W$, such that $ty=(st)(w_{sr})(u_1)$, then
\begin{eqnarray*}
\xi \tilde T_{w_{sr}s}\tilde T_{ty} & = & \xi \tilde T_{w_{sr}s}\tilde T_{(srt)(rw_{sr})(u_1)}\\
                                    & = & \xi^2 \tilde T_{w_{sr}}\tilde T_{(t)(rw_{sr})(u_1)}+\xi \tilde T_{w_{sr}r}\tilde T_{(t)(rw_{sr})(u_1)}
\end{eqnarray*}
Since $w_{sr}trw_{sr}=(w_{sr})(t)(rw_{sr})$, $w_{sr}rtrw_{sr}=(w_{sr}r)(t)(rw_{sr})$,
we have
\begin{eqnarray*}
\xi^2 \tilde T_{w_{sr}}\tilde T_{(t)(rw_{sr})(u_1)}+\xi \tilde T_{w_{sr}r}\tilde T_{(t)(rw_{sr})(u_1)} \\
=\xi^2 \tilde T_{(w_{sr})(t)(rw_{sr})(u_1)}+\xi \tilde T_{(w_{sr}r)(t)(rw_{sr})(u_1)}
\end{eqnarray*}
Meanwhile, since $ty=sty_1\cdots y_{i-1}y_{i+1}\cdots y_k$ and $ty=(st)(w_{sr})(u_1)$,
 we have $$ty_1\cdots y_{i-1}y_{i+1}\cdots y_k=(t)(w_{sr})(u_1)$$
Hence
\begin{eqnarray*}
\tilde T_{w_{sr}s}\tilde T_{ty_1\cdots y_{i-1}y_{i+1}\cdots y_k} & = & \tilde T_{w_{sr}s}\tilde T_{(t)(w_{sr})(u_1)}\\
                                                               & = & \xi \tilde T_{(w_{sr}s)(t)(rw_{sr})(u_1)}+\tilde T_{(w_{sr}sr)(t)(rw_{sr})(u_1)}
\end{eqnarray*}

If there is no $u_1\in W$, such that $ty=(st)(w_{sr})(u_1)$, i.e, $ty=(st)(u_2)$ and $L(u_2)=\{s\}$, then $ty_1\cdots y_{i-1}y_{i+1}\cdots y_k=(t)(u_2)$.
$$\xi \tilde T_{w_{sr}s}\tilde T_{ty}+ \tilde T_{w_{sr}s}\tilde T_{ty_1\cdots y_{i-1}y_{i+1}\cdots y_k}  =
  \xi \tilde T_{w_{sr}sty}+ \tilde T_{w_{sr}s}\tilde T_{tu_2}$$
  By the assumption $L(u_2)=\{s\}$ and $s\notin L(tu_2)$, it is easy too see that $w_{sr}stu_2=(w_{sr}s)(t)(u_2)$.
Hence $$\tilde T_{w_{sr}s}\tilde T_{tu_2}=\tilde T_{w_{sr}stu_2}$$

By Lemma 3.5, in case 2, deg $\tt_{xtr}\tt_y\leq 2$ .

 Case 3: When $x=(x')(w_{st}t)$, then $xtr=(x')(w_{st})(r)$. If $l(xtry)=l(x)+2+l(y)$, then nothing needs to prove.
\begin{eqnarray*}
\tilde T_{xtr}\tilde T_y & = &  \tilde T_{x'w_{st}r}\tilde T_y \\
                         & = &  \tilde T_{x'w_{st}}\tilde T_{ry}
\end{eqnarray*}

We calculate this in the following two conditions:

Condition 1: $R(x')=L(ry)=\{r\}$, here $R(x')=\{r\}$ and $t\notin L(ry)$, since $l(try)=l(y)+2$.
By Lemma 3.6, we have deg $f_{xtr,y,z}\leq 1$.

Condition 2: $L(r,y)=\{s, r\}$.
By 1.1.(d), there exists $u_1\in W$, such that $ry=(w_{sr})(u_1)$.
\begin{eqnarray*}
\tilde T_{x'w_{st}}\tilde T_{ry} & = &  \tilde T_{x'w_{st}}\tilde T_{w_{sr}u_1} \\
                         & = &  \xi\tilde T_{x'w_{st}}\tilde T_{sw_{sr}u_1}+\tilde T_{x'w_{st}s}\tilde T_{sw_{sr}u_1}
\end{eqnarray*}

Since then $R(x')=L(sw_{sr}u_1)=\{r\}$, by Lemma 3.6, we have \\
deg $\xi\tilde T_{x'w_{st}}\tilde T_{sw_{sr}u_1}\leq 2$.

As to the part $\tilde T_{x'w_{st}s}\tilde T_{sw_{sr}u_1}$, $s\notin R(x'w_{st}s)$, then $R(x'w_{st}s)=\{t\}$ or $\{r,t\}$.
In the first case, it is easy to check that $x'stsw_{sr}=(x')(st)(sw_{sr})$, then by Lemma 3.5, $x'stsw_{sr}u_1=(x')(st)(sw_{sr})(u_1)$.
Hence $$\tilde T_{x'w_{st}s}\tilde T_{sw_{sr}u_1}=\tilde T_{x'stsw_{sr}u_1}$$
In the second case, by 1.1.(d), there exists $x''\in W$, such that $x's=(x'')(w_{sr})$.
Then
\begin{eqnarray*}
\tilde T_{x'w_{st}s}\tilde T_{sw_{sr}u_1} & = &\tilde T_{x''w_{sr}t}\tilde T_{sw_{sr}u_1}\\
                                          & = &\xi\tilde T_{x''w_{sr}t}\tilde T_{rsw_{sr}u_1}+\tilde T_{x''w_{sr}rt}\tilde T_{rsw_{sr}u_1}\\
                                          & = &\xi\tilde T_{x''w_{sr}trsw_{sr}u_1}+\tilde T_{x''w_{sr}rtrsw_{sr}u_1}
\end{eqnarray*}

The last equality follows from Lemma 3.5.

Hence in Case 3, we have deg $f_{xtr,y,z}\leq 2$, for all$z\in W$.

 Case 4: When $x=(x')(w_{sr}sr)$, then $xtr=(x')(w_{sr}s)(t)$.
If $l(xtry)=l(x)+2+l(y)$, then nothing needs to prove.
We first calculate $\tilde T_{w_{sr}st}\tilde T_y $.
 Assume that there is an $i<k$, which is minimal, such that $l(w_{sr}sty_1\cdots y_i)<l(w_{sr}sty_1\cdots y_{i-1})$.
By strong exchange condition, and $l(try)=2+l(y)$, we get $srsty_1\cdots y_{i-1}=ty_1\cdots y_{i-1}y_i$, or $rsrty_1\cdots y_{i-1}=ty_1\cdots y_{i-1}y_i$.
$l(ty_1\cdots y_{i-1}y_i)=i+1$. Since $l(srty_1\cdots y_{i-1})=i+2$, then $r\in L(srty_1\cdots y_{i-1})$, by 1.1.(d), there exists $u_1\in W$,s.t
$srty_1\cdots y_{i-1}=(w_{sr})(u_1)$, then $(t)(y_1\cdots y_{i-1})=(rsw_{sr})(u_1)$, which contradicts Corollary 3.2.Hence $rsrty_1\cdots y_{i-1}=ty_1\cdots y_{i-1}y_i$ is impossible.

Then we get $srsty_1\cdots y_{i-1}=ty_1\cdots y_{i-1}y_i$. By the proof of Case 1 in this Lemma,
we see that in fact $ty_1\cdots y_i=(w_{st})(sw_{sr})$. By Lemma 3.5,
\begin{eqnarray*}
\tilde T_{x'w_{sr}st}\tilde T_y & = &  \tilde T_{x'w_{sr}s}\tilde T_{ty_1\cdots y_{i}}\tilde T_{y_{i+1}\cdots y_k} \\
                                & = &  \tilde T_{x'w_{sr}}\tilde T_{(t)(w_{sr})}\tilde T_{y_{i+1}\cdots y_k} \\
                                & = &  \xi \tilde T_{x'w_{sr}}\tilde T_{trw_{sr}y_{i+1}\cdots y_k}+
                                \tilde T_{x'w_{sr}r}\tilde T_{trw_{sr}y_{i+1}\cdots y_k}  \\
                                & = &  \xi \tilde T_{x'w_{sr}rsty}+ \tilde T_{(x')(w_{sr}r)(t)(rw_{sr})(y_{i+1}\cdots y_k)}
\end{eqnarray*}

Hence in Case 4, deg $\tt_{xtr}\tt_y=1$.

In a word, deg $f_{xtr,y,z}\leq 2$, for all $z\in W$.

\medskip

\noindent{\bf  Corollary 3.8.} Let $x,y\in W$. Assume that $R(x)=\{r\}$, $L(y)=\{t\}$, then deg $f_{xw_{st},w_{sr}y,z}\le 2$ for
all $z$ in $W$.

Proof. By the proof of Condition 2 in Case 3 of Lemma 3.7.

\medskip

\noindent{\bf  Lemma 3.9.}
Let $x,y\in W$. Assume that $R(x)=\{s\}$, $L(y)=\{t\}$, then deg $f_{xtr,w_{sr}y,z}\le 3$ for
all $z$ in $W$.

Proof.  $\tilde T_{xtr}\tilde T_{w_{sr}y } = \xi \tilde T_{xtr}\tilde T_{r(w_{sr})(y)}+ \tilde T_{xt}\tilde T_{(rw_{sr})(y)}$.
Obviously, $r\notin R(xt)$.
Since $R(x)=L(rw_{sr}y)=\{s\}$, then deg $\xi \tilde T_{xtr}\tilde T_{r(w_{sr})(y)}\leq 3$, by Lemma 3.7.

Next consider the part $\tilde T_{xt}\tilde T_{(rw_{sr})(y)}$. We have $r\notin R(xt)$, since $R(x)=\{s\}$.
If $R(xt)=\{t\}$, then $\tilde T_{xt}\tilde T_{(rw_{sr})(y)}=\tilde T_{xtrw_{sr}y}$, by lemma 3.5.
If $R(xt)=\{s, t\}$, by 1.1.(d), there exists $x'\in W$, such that $xt=(x')(w_{st})$, then

\begin{eqnarray*}
\tilde T_{xt}\tilde T_{(rw_{sr})(y)} & = & \tilde T_{x'w_{st}}\tilde T_{rw_{sr}y}\\
                                     & = & \xi\tilde T_{x'w_{st}}\tilde T_{srw_{sr}y}+\tilde T_{x'st}\tilde T_{srw_{sr}y}
\end{eqnarray*}
By Lemma 3.6, deg $\xi \tilde T_{x'w_{st}}\tilde T_{(srw_{sr})(y)}\le 2$.
As for the part $\tilde T_{x'st}\tilde T_{srw_{sr}y}$, $s\notin R(x'st)$, since $x'sts=(x')(sts)$, then there are two possibilities.

When $R(x'st)=\{t\}$, by Lemma 3.5, $\tilde T_{x'st}\tilde T_{srw_{sr}y}=\tilde T_{x'stsrw_{sr}y}$.

When $R(x'st)=\{r,t\}$, by 1.1.(d), there exists $x''\in W$, such that $x'st=(x'')(w_{sr})(t)$.
\begin{eqnarray*}
\tilde T_{x'st}\tilde T_{(srw_{sr})(y)} & = & \tilde T_{x''w_{sr}t}\tilde T_{srw_{sr}y}\\
                                        & = & \xi\tilde T_{x''w_{sr}rt}\tilde T_{srw_{sr}y}+\tilde T_{x''w_{sr}r}\tilde T_{trsrw_{sr}y}
\end{eqnarray*}
 Since $R(x''w_{sr}rt)=\{t\}$, by Lemma 3.5,
 $$\xi\tilde T_{x''w_{sr}rt}\tilde T_{srw_{sr}y}=\xi\tilde T_{x''w_{sr}rtsrw_{sr}y}$$
It is easy to see $r\notin L((t)(rsrw_{sr})(y))$, otherwise it contradicts to the fact $L(y)=\{t\}$.
When $L(trsrw_{sr}y)=\{t\}$, by Lemma 3.5,
$$\tilde T_{x''w_{sr}r}\tilde T_{trsrw_{sr}y}=\tilde T_{x''w_{sr}rtrsrw_{sr}y}$$
When $L(trsrw_{sr}y)=\{s,t\}$, meanwhile $m_{sr}=7$ and there exists $y'\in W$, such that $y=(t)(rw_{sr})(y')$.
Then $trsrw_{sr}y'=tstrstsw_{sr}y''$,
\begin{eqnarray*}
\tilde T_{x''w_{sr}r}\tilde T_{(trsrw_{sr})(y)} & = & \tilde T_{x''w_{sr}r}\tilde T_{stsrstsw_{sr}y''}\\
                                        & = & \xi\tilde T_{x''w_{sr}r}\tilde T_{tsrstw_{sr}y''}+\tilde T_{x''w_{sr}rs}\tilde T_{tsrstw_{sr}y''}
\end{eqnarray*}
Since $L(tsrstw_{sr}y'')=\{t\}$, by Lemma 3.5,
$$\xi\tilde T_{x''w_{sr}r}\tilde T_{tsrstw_{sr}y''}+\tilde T_{x''w_{sr}rs}\tilde T_{tsrstw_{sr}y''}=\xi\tilde T_{x''w_{sr}rtsrstw_{sr}y''}+\tilde T_{x''w_{sr}rstsrstw_{sr}y''}$$
In a word, deg $f_{xtr,w_{sr}y,z}\le 3$, for all $z$ in $W$.

\medskip

\noindent{\bf  Lemma 3.10.}
Let $x,y\in W$. Assume that $R(x)=\{s\}$, $L(y)=\{r\}$, then deg $f_{xtr,w_{st}y,z}\le 4$ for
all $z$ in $W$.

Proof.
$\tilde T_{xtr}\tilde T_{(tst)(y)} = \xi\tilde T_{xtr}\tilde T_{sty}+\tilde T_{xr}\tilde T_{sty}$.
Obviously, $t,r\notin L(sty)$.
$R(x)=L(sty)=\{s\}$, deg $\xi\tilde T_{xtr}\tilde T_{sty}\leq 3$, by Lemma 3.7.
As to part $\tilde T_{xr}\tilde T_{sty}$, since $L(y)=\{r\}$, write $y=ry_1$, $L(y_1)=\{s\}$.
$\tilde T_{xr}\tilde T_{sty}=\tilde T_{xr}\tilde T_{stry_1}$.

When $R(xr)=\{r\}$, it is easy to check that $t\notin R(xrs)$.

1) $R(xrs)=\{s\}$, by Lemma 3.7, deg $\tilde T_{xr}\tilde T_{sty}\leq 2$.

2) $R(xrs)=\{s,r\}$, by Lemma 3.9, deg $\tilde T_{xr}\tilde T_{sty}\leq 3$

When $R(xr)=\{s,r\}$, there exists $x'\in W$, such that $xr=(x')(w_{sr})$,
$$\tilde T_{xr}\tilde T_{sty}=\tilde T_{x'w_{sr}}\tilde T_{stry_1}=\xi\tilde T_{x'w_{sr}}\tilde T_{try_1}+\tilde T_{x'w_{sr}s}\tilde T_{try_1}$$
By Lemma 3.9, deg $\xi\tilde T_{x'w_{sr}}\tilde T_{try_1}\leq 4$.
$$\tilde T_{x'w_{sr}s}\tilde T_{try_1}=\xi\tilde T_{x'w_{sr}sr}\tilde T_{try_1}+\tilde T_{x'w_{sr}sr}\tilde T_{ty_1}$$
Since $R(x'w_{sr}sr)=\{s\}=L(y_1)$, deg $\xi\tilde T_{x'w_{sr}sr}\tilde T_{try_1}\leq 3$, by Lemma 3.7.

Finally we consider the part $\tilde T_{x'w_{sr}sr}\tilde T_{ty_1}$.
Obviously, $r\notin L(ty_1)$.

If $L(ty_1)=\{t\}$, by Lemma 3.5, $\tilde T_{x'w_{sr}sr}\tilde T_{ty_1}=\tilde T_{x'w_{sr}srty_1}$.

If $L(ty_1)=\{s,t\}$, by 1.1.(d), there exists $y_2\in W$, $L(y_2)=\{r\}$, such that $ty_1=(w_{st})(y_2)$.
$$\tilde T_{x'w_{sr}sr}\tilde T_{ty_1}=\tilde T_{x'w_{sr}sr}\tilde T_{w_{st}y_2}=\xi \tilde T_{x'w_{sr}sr}\tilde T_{tsy_2}+\tilde T_{x'w_{sr}srs}\tilde T_{tsy_2}$$
Obviously $s\notin L(tsy_2)$.
If $L(tsy_2)=\{t\}$, then it is easy to check that $L(sy_2)=\{s\}$.
By Lemma 3.5, $\xi \tilde T_{x'w_{sr}sr}\tilde T_{tsy_2}=\xi \tilde T_{x'w_{sr}srtsy_2}$.
Since $R(x'w_{sr}srsr)=\{s\}$, $\tilde T_{x'w_{sr}srs}\tilde T_{tsy_2}=\tilde T_{x'w_{sr}srsr}\tilde T_{trsy_2}$ then by Lemma 3.7, deg $\tilde T_{x'w_{sr}srs}\tilde T_{tsy_2}\leq 2$.
If $L(tsy_2)=\{t,r\}$, by 1.1.(d), there exists $y_3\in W$, such that $sy_2=(w_{sr})(y_3)$.
Since $\xi \tilde T_{x'w_{sr}sr}\tilde T_{tsy_2}=\xi \tilde T_{x'w_{sr}srt}\tilde T_{w_{sr}y_3}$, and $R(x'w_{sr}srt)=\{t\}$, by Lemma 3.5, $\xi \tilde T_{x'w_{sr}sr}\tilde T_{tsy_2}=\xi \tilde T_{x'w_{sr}srtsy_2}$.
Since $\tilde T_{x'w_{sr}srs}\tilde T_{tsy_2}= \tilde T_{x'w_{sr}srs}\tilde T_{tw_{sr}y_3}=\tilde T_{(x'w_{sr}srsr)(rt)}\tilde T_{w_{sr}y_3}$,
by Lemma 3.9, deg $\tilde T_{x'w_{sr}srs}\tilde T_{tsy_2}\leq 3$.

Hence we can conclude that  deg $f_{xtr,w_{st}y,z}\le 4$ for
all $z$ in $W$.

\medskip

Let $P$ be the parabolic subgroup of W generated by $s$ and $r$.

\medskip

\noindent{\bf  Lemma 3.11.} Assume that $w$, $u$ are elements of $P$. Then deg $f_{w,u,v}\le l(v)$ for $v\in P$ and deg $f_{w,u,v}=0$ if $v\notin P$.

Proof. Refer to [X].

\medskip

\noindent{\bf  Lemma 3.12.} Let $x,y\in W$. Let $x_1$ (resp. $y_1$)be the element in the coset $xP$ (resp. $Py$) with minimal length. Let $w,u\in P$ be such that $x=x_1w$, $y=uy_1$. When $l(w),l(u)\geq 1$ and $l(w)+l(u)\geq 3$,
then deg $f_{x,y,z}\le m_{sr}$ for
all $z$ in $W$.

Proof. We use induction on min $\{l(x),l(y)\}$. When min $\{l(x),l(y)\}\leq m_{sr}$,
the lemma is clear. Next assume that $k> m_{sr}$.
By the assumption, we have
$$\tt_x\tt_y=\sum_{v\in P}f_{w,u,v}\tt_{x_1v}\tt_{y_1}.$$
By Lemma 3.11, deg$f_{w,u,v}\le l(v)$ and $v\in P$ if $f_{w,u,v}\ne
0$. If $l(v)\ge 5$, by Lemma 3.5, $l(x_1vy_1)=l(x_1v)+l(y_1)$. Hence
$\tt_{x_1w}\tt_{y_1}=\tt_{x_1wy_1}$.

If  $l(v)=0$

$R(x_1)=L(y_1)=\{t\}$. Write $x_1=(x_2)(t)$, $y_1=(t)(y_2)$, here $R(x_2)=L(y_2)=\{s\}$.
$$\tt_{x_1}\tt_{y_1}=\tt_{x_2t}\tt_{ty_2}=\xi\tt_{x_2t}\tt_{y_2}+\tt_{x_2}\tt_{y_2}$$
Write $x_2=x_3s$, $y_2=sy_3$, then it is easy to check that $R(x_3)=L(y_3)=\{r\}$ .
Hence, by Lemma 3.6, $$\xi\tt_{x_2t}\tt_{y_2}=\xi\tt_{x_3sts}\tt_{y_3}$$  deg $\xi\tt_{x_3sts}\tt_{y_3}\leq 2$.
$\tt_{x_2}\tt_{y_2}=\tt_{x_3s}\tt_{sy_3}$, by induction hypotheses, deg $\tt_{x_2}\tt_{y_2}\leq m_{sr}$.
Hence deg $\tt_{x_1}\tt_{y_1}\leq m_{sr}$.

Write $x_1=x_2rst$, $y_1=tsry_2$, since $R(x_1)=L(y_1)=\{t\}$.
It is to check that $R(x_2rs)=L(sry_2)=\{s\}$, $R(x_2r)=L(ry_2)=\{r\}$, $R(x_2)=L(y_2)=\{s\}$.

If $l(v)=1$

1) $v=r$.
$$\tt_{x_1r}\tt_{y_1}=\tt_{x_2rstr}\tt_{tsry_2}=\xi\tt_{x_2rstr}\tt_{sry_2}+\tt_{x_2rsr}\tt_{sry_2}$$
By lemma 3.7, deg $\xi\tt_{x_2rstr}\tt_{sry_2}\leq 3$.
Since $\xi\tt_{x_2rsr}\tt_{sry_2}=\tt_{x_2rsr}\tt_{rsry_2}-\tt_{x_2rs}\tt_{sry_2}$.Here $l(x_1)=l(x_2)+3$, $l(y_1)=l(y_2)+3$,
and $l(x_1)\leq l(x)-1$, $l(y_1)\leq l(y)-1$, hence we can use induction hypotheses to $\xi\tt_{x_2rsr}\tt_{sry_2}$, and the lemma is true then.

2) $v=s$. It is easy to check that $R(x_1s)=\{s,t\}$, by 1.1.(d), write $x_1s=(x_3)(w_{st})$, $x_3=x_2r$, $R(x_3)=\{r\}$.
\begin{eqnarray*}
\tt_{x_1s}\tt_{y_1} & = & \tt_{x_3w_{st}}\tt_{tsry_2}\\
                    & = & \xi\tt_{x_3tst}\tt_{sry_2}+\tt_{x_3ts}\tt_{sry_2}\\
                    & = & \xi^2\tt_{x_3sts}\tt_{ry_2}+\xi\tt_{x_3st}\tt_{ry_2}+\xi\tt_{x_3ts}\tt_{ry_2}+\tt_{x_3t}\tt_{ry_2}
\end{eqnarray*}
Since $R(x_3)=L(ry_2)=\{r\}$, by Lemma 3.6, deg $\xi^2\tt_{x_3sts}\tt_{ry_2}\leq 3$.
Since $R(x_2rs)=L(sry_2)=\{s\}$, $\xi\tt_{x_3st}\tt_{ry_2}=\xi\tt_{x_2rst}\tt_{ry_2}=\xi\tt_{x_2rstr}\tt_{y_2}$.
Since $R(x_2rs)=L(y_2)=\{s\}$, by Lemma 3.7, deg $\xi\tt_{x_3st}\tt_{ry_2}\leq 3$.
Since $\xi\tt_{x_3ts}\tt_{ry_2}=\xi\tt_{x_2rt}\tt_{sry_2}$, $R(x_2)=L(sry_2)=\{s\}$, then by Lemma 3.7, deg $\xi\tt_{x_3ts}\tt_{ry_2}\leq 3$.

As for the part $\tt_{x_3t}\tt_{ry_2}$, $\xi \tt_{x_3t}\tt_{ry_2}=\xi\tt_{x_2rt}\tt_{ry_2}=\tt_{x_2rt}\tt_{try_2}-\tt_{x_2r}\tt_{ry_2}$
Apply induction hypotheses, we see that deg $\tt_{x_2r}\tt_{ry_2}\leq m_{sr}$.
Then we deal with the following, $$\tt_{x_2rt}\tt_{ry_2}=\xi\tt_{x_2rt}\tt_{y_2}+\tt_{x_2t}\tt_{y_2}$$
By Lemma 3.7, we have deg $\xi^2\tt_{x_2rt}\tt_{y_2}\leq 4$.

Finally, claim that deg $\tt_{x_2t}\tt_{y_2}\leq m_{sr}-1$, i.e, deg $\xi \tt_{x_2t}\tt_{y_2}\leq m_{sr}$, when $R(x_2)=L(y_2)=\{s\}$.
(\textbf{Notice} The claim here will be used in the proof of Lemma 3.13. )
Choose suitable $x'$, $y'\in W$, $R(x')=L(y')=\{r\}$.

1) $x_2t=(x')(w_{st})$, $y_2=sty'$
\begin{eqnarray*}
\xi \tt_{x_2t}\tt_{y_2} & = & \xi \tt_{x'w_{st}}\tt_{sty'}\\
                        & = & \xi^2 \tt_{x'w_{st}}\tt_{ty'}+\xi \tt_{x'st}\tt_{ty'}\\
                        & = & \xi^3 \tt_{x'w_{st}}\tt_{y'}+\xi^2 \tt_{x'ts}\tt_{y'}+\xi ^2\tt_{x'st}\tt_{y'}+\xi \tt_{x's}\tt_{y'}
 \end{eqnarray*}
 By Lemma 3.6, deg $\xi^3 \tt_{x'w_{st}}\tt_{y'}\leq 4$.
 $\xi^2 \tt_{x'ts}\tt_{y'}=\xi^2 \tt_{x't}\tt_{sy'}$,since $t\notin L(sy')$,$s\in L{sy'}$, by Lemma 3.6, or Lemma 3.9, deg $\xi^2 \tt_{x'ts}\tt_{y'}\leq 5$.
 As the same reason, deg $\xi^2 \tt_{x'st}\tt_{y'}\leq 5$.
 then apply induction hypotheses to the left part $\xi \tt_{x's}\tt_{y'}$, which is equal to $\tt_{x's}\tt_{sy'}-\tt_{x'}\tt_{y'}$.

 2) $x_2t=(x')(w_{st})$, $y_2=sy'$
\begin{eqnarray*}
\xi \tt_{x_2t}\tt_{y_2} & = & \xi \tt_{x'w_{st}}\tt_{sy'}\\
                        & = & \xi^2 \tt_{x'w_{st}}\tt_{y'}+\xi \tt_{x'st}\tt_{y'}\\
                        & = & \xi^2 \tt_{x'w_{st}}\tt_{y'}+\xi \tt_{x't}\tt_{sy'}
 \end{eqnarray*}
 By Lemma 3.6, deg $\xi^2 \tt_{x'w_{st}}\tt_{y'}\leq 3$.
 By Lemma 3.7, deg $\xi \tt_{x't}\tt_{sy'}\leq 2$.

 3) $x_2t=(x')(st)$, $y_2=sty'$
 \begin{eqnarray*}
\xi \tt_{x_2t}\tt_{y_2} & = & \xi \tt_{x'st}\tt_{sty'}\\
                        & = & \xi \tt_{x'w_{st}}\tt_{ty'}\\
                        & = & \xi^2 \tt_{x'w_{st}}\tt_{y'}+\xi \tt_{x'ts}\tt_{y'}
 \end{eqnarray*}
 By Lemma 3.6,  deg $\xi^2 \tt_{x'w_{st}}\tt_{y'}\leq 3$.
 $\xi \tt_{x'ts}\tt_{y'}=\xi \tt_{x't}\tt_{sy'}$, $t\notin L(sy')$, $s\in L(sy')$,
 by Lemma 3.7 or Lemma 3.9 deg $\xi \tt_{x'ts}\tt_{y'}\leq 4$.

 4) $x_2t=(x')(st)$, $y_2=sy'$
 $$\xi \tt_{x_2t}\tt_{y_2} = \xi \tt_{x'st}\tt_{sy'}=\xi \tt_{x'sts}\tt_{y'}$$
 By Lemma 3.6, deg $\xi \tt_{x_2t}\tt_{y_2}\leq 1 $.

 Hence deg $\xi \tt_{x_2t}\tt_{y_2}\leq m_{sr}$.

If $l(v)=2$, $v=sr$ or $v=rs$.
$$\tt_{x_1sr}\tt_{y_1}=\tt_{x_2rsts}\tt_{rtsry_2}$$
By Lemma 3.10, deg $\tt_{x_1sr}\tt_{y_1}\leq 4$.
As the same, deg $\tt_{x_1rs}\tt_{y_1}\leq 4$.

If $l(v)=3$

1) $v=srs$
\begin{eqnarray*}
\tt_{x_1srs}\tt_{y_1} & = & \tt_{x_2rstsrs}\tt_{tsry_2}\\
                      & = & \tt_{x_2rstsr}\tt_{stsry_2}\\
                      & = & \xi \tt_{x_2rtsr}\tt_{stsry_2}+\tt_{x_2rtsr}\tt_{stry_2}\\
                      & = & \xi \tt_{x_2rt}\tt_{srstsry_2}+\tt_{x_2rt}\tt_{srstry_2}\\
                      & = & \xi \tt_{x_2rtsrsr}\tt_{stsry_2}+\tt_{x_2rt}\tt_{srstry_2}
\end{eqnarray*}
It is easy to check $L(srstsry_2)=L(srstry_2)=\{s\}$, then by Lemma 3.7, deg $\xi \tt_{x_2rt}\tt_{srstsry_2}\leq 3$,
deg $\tt_{x_2rt}\tt_{srstry_2}\leq 2$.
Or since $R(x_2rtsrsr)=L(ry_2)=\{r\}$, by Lemma 3.5,  deg $\xi \tt_{x_2rt}\tt_{srstsry_2}\leq 2$.

2) $v=rsr$
\begin{eqnarray*}
\tt_{x_1rsr}\tt_{y_1} & = & \tt_{x_2rstrsr}\tt_{tsry_2}\\
                      & = & \tt_{x_2rsrts}\tt_{trsry_2}\\
\end{eqnarray*}
When $R(x_2rsrts)=\{s,t\}$, then $R(x_2rsr)=\{s,r\}$, by 1.1.(d), write $x_2rsr=(x_4)(w_{sr})$.
\begin{eqnarray*}
\tt_{x_1rsr}\tt_{y_1} & = & \tt_{x_2rsrts}\tt_{trsry_2}\\
                      & = & \xi \tt_{x_4w_{sr}ts}\tt_{rsry_2}+\tt_{x_4w_{sr}sts}\tt_{rsry_2}
\end{eqnarray*}
Since $R((x_4)(w_{sr}s))=\{r\}$, $L(rsry_2)=\{r\}$ or $\{s,r\}$ , then by Lemma 3.6, or Corollary 3.8,
deg $\xi \tt_{x_4w_{sr}ts}\tt_{rsry_2}\leq 3$.

When $L(rsry_2)=\{r\}$, $\tt_{x_4w_{sr}sts}\tt_{rsry_2}=\tt_{x_4w_{sr}st}\tt_{srsry_2}$, by Lemma 3.6, $R(x_4w_{sr}sr)=\{s\}$,
$L(srsry_2)=\{s\}$, or $\{s,r\}$, then by Lemma 3.7, or 3.9, deg $\tt_{x_4w_{sr}st}\tt_{srsry_2}\leq 3$.

When $L(rsry_2)=\{s,r\}$, write $rsry_2=(w_{sr})(y_3)$, by Lemma 1.3.
\begin{eqnarray*}
\tt_{x_4w_{sr}sts}\tt_{rsry_2} & = & \tt_{(x_4)(w_{sr}s)(ts)}\tt_{w_{sr}y_3}\\
                               & = & \xi\tt_{(x_4)(w_{sr}s)(t)}\tt_{w_{sr}y_3}+\tt_{x_4w_{sr}st}\tt_{sw_{sr}y_3} \\
                               & = & \xi^2\tt_{(x_4)(w_{sr}sr)(t)}\tt_{w_{sr}y_3}+\xi\tt_{(x_4)(w_{sr}sr)(t)}\tt_{(rw_{sr})(y_3)}\\
                               &   & +\xi\tt_{(x_4)(w_{sr}s)(t)}\tt_{(rsw_{sr})(y_3)}+\tt_{(x_4)(w_{sr}sr)(t)}\tt_{(rsw_{sr})(y_3)}
\end{eqnarray*}
By Lemma 3.5, we see that deg $\tt_{x_4w_{sr}sts}\tt_{rsry_2}=2$.

If $l(v)=4$, $v=srsr$ or $rsrs$.
When $v=srsr$,
$$\tt_{x_1srsr}\tt_{y_1} = \tt_{x_2rstsrs}\tt_{trsry_2}$$
It is easy to check that $R(x_2rstsrs)=L(sry_2)=\{s\}$, hence by Lemma 3.7,
deg $\tt_{x_1srsr}\tt_{y_1}\leq 2$.
As the same reason, deg $\tt_{x_1}\tt_{rsrsy_1}\leq 2$.

Hence the lemma is proved.

\medskip

\noindent{\bf Theorem 3.13.}
 $(W,S)$ is a Coxeter group, $S=\{r,s,t\}$, $m_{sr}\geq 7, m_{st}=3,rt=tr$.
Then deg $f_{x,y,z}\le m_{sr}$ for
all $x,y,z$ in $W$.

Proof. $\forall x,y \in W$, we discuss it in the following 6 cases.

1) $R(x)=\{t\}$

When $L(y)=\{t\}$, write $x=(x_0)(st)$, $y=(ts)(y_0)$, here $R(x_0)=L(y_0)=\{r\}$, $R(x_0s)=L(sy_0)=\{s\}$.
$$\tt_x\tt_y=\tt_{x_0st}\tt_{tsy_0}=\xi\tt_{x_0st}\tt_{sy_0}+\tt_{x_0s}\tt_{sy_0}$$
By the notice in the proof of Lemma 3.12, deg $\xi\tt_{x_0st}\tt_{sy_0}\leq m_{sr}$.
By Lemma 3.12, deg $\tt_{x_0s}\tt_{sy_0}\leq m_{sr}$.

When $L(y)=\{s,t\}$, write $y=(w_{st})(y_1)$, $L(y_1)=\{r\}$.
 \begin{eqnarray*}
 \tt_x\tt_y & = & \tt_{x_0st}\tt_{stsy_1}\\
            & = & \xi\tt_{x_0s}\tt_{tsty_1}+\tt_{x_0s}\tt_{sty_1}\\
            & = & \xi^2 \tt_{x_0}\tt_{stsy_1}+\xi\tt_{x_0}\tt_{tsy_1}+\xi\tt_{x_0}\tt_{sty_1}+\tt_{x_0}\tt_{ty_1}
  \end{eqnarray*}
  By Lemma 3.6, deg $\xi^2 \tt_{x_0}\tt_{stsy_1}\leq 3$.
  By Lemma 3.7 or Lemma 3.9, deg $\xi\tt_{x_0}\tt_{tsy_1}\leq 4$, deg $\xi\tt_{x_0}\tt_{sty_1}\leq 4$.
  As for $\tt_{x_0}\tt_{ty_1}$, it will be proved later.

 When $L(y)=\{s,r\}$, by Lemma 3.6, it is done.

 When $L(y)=\{t,r\}$, this will be proved in 2).

 When $L(y)=\{s\}$, this will be proved in 3).

 When $L(y)=\{r\}$, by Lemma 3.7, it is done.

2) $R(x)=\{t,r\}$ ,write $x=(x_2)(tr)$.

 When $L(y)=\{s\}$, by Lemma 3.7,deg $f_{x,y,z}\le 2$, hence deg $f_{x,y,z}\le m_{sr}$.

 When $L(y)=\{s,t\}$, by Lemma 3.10, deg $f_{x,y,z}\le 4$.

 When $L(y)=\{s,r\}$, by Lemma 3.12, this is done.

 When $L(y)=\{r\}$, write $y=(rs)(y_2)$, here $L(sy_2)=\{s\}$.
 Hence by Lemma 3.12, this is done. So is $\tt_{x_0}\tt_{ty_1}$ in 1).

 When $L(y)=\{r,t\}$, write $y=(tr)(y_4)$, $L(y_4)=\{s\}$.
 $$\tt_x\tt_y=\tt_{x_2tr}\tt_{try_4}=\xi^2\tt_{x_2tr}\tt_{y_4}+\xi\tt_{x_2t}\tt_{y_4}+\xi\tt_{x_2r}\tt_{y_4}+\tt_{x_2}\tt_{y_4}$$
 By Lemma 3.7, deg $\xi^2\tt_{x_2tr}\tt_{y_4}\leq 4$.
 By the proof of Lemma 3.12,deg $\xi\tt_{x_2t}\tt_{y_4}\leq m_{sr}$.
 Since $\xi\tt_{x_2r}\tt_{y_4}+\tt_{x_2}\tt_{y_4}=\tt_{x_2r}\tt_{ry_4}$, by Lemma 3.12 deg $\tt_{x_2r}\tt_{ry_4}\leq m_{sr}$.

 When $L(y)=\{t\}$, write $y=(ts)(y_5)$, $L(y_5)=\{r\}$, $L(sy_5)=\{s\}$.
 $$\tt_x\tt_y=\tt_{x_2tr}\tt_{tsy_5}=\xi\tt_{x_2tr}\tt_{sy_5}+\tt_{x_2r}\tt_{sy_5}$$
 By Lemma 3.7, deg $\xi\tt_{x_2tr}\tt_{sy_5}\leq 3$.
 By Lemma 3.12, deg $\tt_{x_2r}\tt_{sy_5}\leq m_{sr}$.

 3) $R(x)=\{s\}$, we deal this in two conditions.

Condition 1: $x=(x_3)(ts)$, $R(x_3)=\{r\}$ .

When $L(y)=\{r\}$, or $y=(sr)(y_6)$, including $L(y)=\{s,r\}$, by Lemma 3.12, they are done.

When $L(y)=\{s,t\}$, and write $y=(w_{st})(y_7)$, $L(y_7)=\{r\}$.
 \begin{eqnarray*}
 \tt_x\tt_y & = & \tt_{x_3ts}\tt_{stsy_7}\\
            & = & \xi\tt_{x_3t}\tt_{stsy_7}+\tt_{x_3t}\tt_{tsy_7}\\
            & = & \xi^2\tt_{x_3}\tt_{stsy_7}+\xi\tt_{x_3}\tt_{sty_7}+\xi\tt_{x_3t}\tt_{sy_7}+\tt_{x_3}\tt_{sy_7}
\end{eqnarray*}
By Lemma 3.6, deg $\xi^2\tt_{x_3}\tt_{stsy_7}\leq 3$.
By Lemma 3.7, or Lemma 3.9, deg $\xi\tt_{x_3}\tt_{sty_7}\leq 4$.
By Lemma 3.7, or Lemma 3.9, deg $\xi\tt_{x_3t}\tt_{sy_7}\leq 4$.
By Lemma 3.12, deg $\tt_{x_3}\tt_{sy_7}\leq m_{sr}$.

When $L(y)=\{s\}$ and $y=(st)(y_8)$, $L(y_8)=\{r\}$
\begin{eqnarray*}
 \tt_x\tt_y & = & \tt_{x_3ts}\tt_{sty_8}\\
            & = & \xi\tt_{x_3t}\tt_{sty_8}+\tt_{x_3t}\tt_{ty_8}\\
            & = & \xi\tt_{x_3}\tt_{tsty_8}+\xi\tt_{x_3}\tt_{ty_8}+\tt_{x_3}\tt_{y_8}
\end{eqnarray*}
By Lemma 3.6, deg $\xi\tt_{x_3}\tt_{tsty_8}\leq 2$.
$\xi\tt_{x_3}\tt_{ty_8}=\tt_{x_3t}\tt_{ty_8}-\tt_{x_3}\tt_{y_8}$
Since $R(x_3t)=L(ty_8)=\{t,r\}$, by what we have proved before, its degree is less than $m_{sr}$.
Since $R(x_3)=L(y_8)=\{r\}$, by Lemma 3.12, this is done.

When $L(y)=\{t\}$, write $y=(ts)(y_9)$, $L(y_9)=\{r\}$.
\begin{eqnarray*}
 \tt_x\tt_y & = & \tt_{x_3ts}\tt_{tsy_9}\\
            & = & \tt_{x_3tst}\tt_{sy_9}\\
            & = & \xi\tt_{x_3sts}\tt_{y_9}+\tt_{x_3s}\tt_{ty_9}
\end{eqnarray*}
By Lemma 3.6, deg $\xi\tt_{x_3sts}\tt_{y_9}\leq 2$
By Lemma 3.7, or Lemma 3.9, deg $\tt_{x_3s}\tt_{ty_9}\leq 3$.

When $l(y)=\{t,r\}$, which has been already done in 2).

Condition 2: When $x=(x_4)(srs)$, $R(x_4sr)=\{r\}$
It can be dealt with Lemma 3.12. Hence it is done.

4) $R(x)=\{r\}$

It is easy to check that by Lemma 3.12,  $L(y)=\{r\}$, $\{s,r\}$, $\{s,t\}$, $\{r,t\}$,
and $\{s\}$ are done.

When $L(y)=\{t\}$, it is done in 1).

5) $R(x)=\{s,r\}$

 For all $y\in W$, this is done by Lemma 3.12.

6) $R(x)=\{s,t\}$

It is easy to check that by Lemma 3.12,  $L(y)=\{r\}$ and $\{s,r\}$ are done.

When $L(y)=\{t\}$, it is done in 1).

When $L(y)=\{r,t\}$, it is done in 2).

When $L(y)=\{s\}$, it is done in 3).

When $L(y)=\{s,t\}$.
Write $x=(x_5)(w_{st})$, $y=(w_{st})(y_{10})$.
\begin{eqnarray*}
 \tt_x\tt_y & = & \tt_{x_5sts}\tt_{stsy_{10}}\\
            & = & \xi^3\tt_{x_5sts}\tt_{y_{10}}+\xi^2\tt_{x_5ts}\tt_{y_{10}}+\xi^2\tt_{x_5st}\tt_{y_{10}} \\
            &   & +\xi\tt_{x_5s}\tt_{y_{10}}+\xi\tt_{x_5t}\tt_{y_{10}}+\xi\tt_{x_5sts}\tt_{y_{10}}+\tt_{x_5}\tt_{y_{10}}
\end{eqnarray*}
By Lemma 3.6, deg $ \xi^3\tt_{x_5sts}\tt_{y_{10}}\leq 4$, deg $\xi\tt_{x_5sts}\tt_{y_{10}}\leq 2$.
By Lemma 3.7, or Lemma 3.9, deg $\xi^2\tt_{x_5ts}\tt_{y_{10}}\leq 5$, deg $\xi^2\tt_{x_5st}\tt_{y_{10}}\leq 5$.
By Lemma 3.12, deg $\tt_{x_5}\tt_{y_{10}}\leq m_{sr}$.

 Hence the theorem is proved.

\medskip

\section{The case $m_{sr}\geq 5$ and  $m_{st}\geq 4 $}

In this section $(W,S)$ is a Coxeter group of rank 3, $S=\{s,t,r\}$, $rt=tr$.
Firstly, we assume that $m_{sr}\geq 4$ and $m_{st}\geq 4 $.

\medskip

\noindent{\bf Lemma 4.1.} Keep the assumptions and  notations
above. There is no element $w$ in $W$ such that $w=(w_1)(r)=(w_2)(ts)$ .

Proof. We use induction on $l(w)$. When $l(w)=0,1,2,3$, the lemma is
clear. Now assume that the lemma is true for $u$ with
$l(u)\leq l(w)-1$. Since $r,s\in R(w)$. By 1.1.(d),
$w=(w_3)(w_{sr})$ for some $w_3\in W$.
So we get $w_1=(w_3)(w_{sr}r),w_2t=(w_3)(w_{sr}s)$. Then $r,t\in R(w_2t)$.
By 1.1.(d),
$w_2t=w_3w_{sr}s=(w_4)(tr)$ for some $w_4\in W$.
$w_2=w_4r$, $(\tilde{w_3})(rs)=(w_4)(t)$ for some $\tilde{w_3}\in W$, since $m_{sr}\geq 4$.
By calculation, there exists  $w_5\in W$, such that $(\tilde{w_3})(rs)=(w_4)(t)=(w_5)(w_{st})$, by Lemma 1.3.
That is $(\tilde{w_3})(r)=\tilde{w_5}(tst)$, here $\tilde{w_5}tst=w_5w_{st}s$.
Then there exists $w_6\in W$, such that $(\tilde{w_5})(tst)=(w_6)(tr)$, by Lemma 1.3.
Hence $(\tilde{w_5})(ts)=(w_6)(r)$, which by induction hypothesis is impossible. The
lemma is proved.

\medskip

\noindent{\bf Corollary 4.2.} There is no element $w$ in $W$ such that $w=(w_1)(t)=(w_2)(rs)$.

Proof. From the proof of Lemma 4.1.

\medskip
\noindent{\bf Lemma 4.3.} There is no element $w$ in $W$ such that

\noindent (a) $w=(w_1)(r)=(w_2)(sts)$.

\noindent (b)  $w=(w_1)(r)=(w_2)(tst)$.

\noindent (c) $w=(w_1)(t)=(w_2)(srs)$.

\noindent (d) $w=(w_1)(t)=(w_2)(rsr)$.

Proof. We only have to deal with (a) and (b).

 By Lemma 4.1, (a) is done.

 We use induction on $l(w)$. When $l(w)=0,1,2,3$, the lemma is
clear. Now assume that the lemma is true for $u$ with
$l(u)\leq l(w)-1$. Since $r,t\in R(w)$. By 1.1.(d),
$w=(w_3)(tr)$ for some $w_3\in W$.
So we get $w_1=(w_3)(t),(w_2)(ts)=(w_3)(r)$, which contradicts Lemma 4.1.
Hence (b) is proved.

\medskip

\noindent{\bf Lemma 4.4.}
There is no element $w$ in $W$ such that
 $w=(w_1)(sr)=(w_2)(st)$.

 Proof. We use induction on $l(w)$. When $l(w)=0,1,2,3$, the lemma is
clear. Now assume that the lemma is true for $u$ with
$l(u) \leq l(w)-1$. Since $r,t\in R(w)$. By 1.1.(d),
$w=(w_3)(tr)$ for some $w_3\in W$.
So we get $(w_1)(s)=(w_3)(t),(w_2)(s)=(w_3)(r)$. Then  by 1.1.(d),
$w_3=(w_4)(w_{st}t)$ for some $w_4\in W$.
$w_3=(w_5)(w_{sr}r)$ for some $w_4\in W$. Hence $(\tilde{w_4})(sts)=(\tilde{w_5})(srs)$,  since $m_{sr}, m_{st}\geq 4$.Here $(\tilde{w_4})(sts)=(w_4)(w_{st}t)$,  $(\tilde{w_5})(srs)=(w_5)(w_{sr}r)$. By induction hypothesis, $(\tilde{w_4})(st)=(\tilde{w_5})(sr)$ is
impossible , the lemma is proved.

\medskip

\textbf{The Notation}, let $\{\alpha,\beta\}=\{t,r\}$.

\medskip
\noindent{\bf Lemma 4.5.}
Let $x,y$ be elements in $W$, and $w$ be an element in the parabolic subgroup
 generated by the two simple reflections $s,\alpha$, $l(w)\ge 4$ and $s,\alpha$ are not in $R(x)\cup L(y)$. Then

 \noindent (a)$l(xwy)=l(x)+l(w)+l(y)$.

 \noindent (b)$R(xwy)=R(wy)$.

 \noindent (c)$L(xwy)=L(xw)$.

Proof. It is clear that $xw=(x)(w)$, and $wy=(w)(y)$. Note that (b) and (c) are equivalent.
We use induction on $l(y)$ to prove (a) and (b).

 When $l(y)=0$, since $l(w)\geq 4$, by Lemma 4.3, $\beta\notin R(xw)$. When $R(w)=\{\alpha,s\}$, $R(xw)=\{\alpha,s\}$. When $R(w)=\{s\}$ or $\{\alpha\}$, since $R(x)=\{\beta\}$, $R(w)=R(xw)$.
When $l(y)=1$, i.e., $y=\beta$. If $\alpha\in R(w)$, then $R(xw\beta)=R(w\beta)=\{t,r\}$.
If $R(w)=\{s\}$, then it is easy to check that $R(xw\beta)=R(w\beta)=\{\beta\}$.
when $l(y)=2$, i.e., $y=\beta s$. By Lemma 4.1, Corollary 4.2, $\alpha\notin R(xwy)$. If $\beta\in R(xw\beta s)$,
then $s\in R(xw)=R(w)$, since $l(w_{s\beta})\geq 4$, $\beta\in R(\tilde{x}\alpha s \alpha)$, here $(\tilde{x})(\alpha s \alpha)=(x)(ws)$, which contradicts Lemma 4.3. Hence $R(xw\beta s)=R(w\beta s)=\{s\}$.
Next assume that $l(y)\geq 3$.
Assume that the lemma is true when $l(y)\leq k-1$, $k\geq 3$. When $l(y)=k$, Write $y=y_1\cdots y_k$, reduced decomposition.
The induction hypothesis says that $R(xwy_1\cdots y_i)=R(wy_1\cdots y_i)$ and $l(xwy_1\cdots y_i)=l(x)+l(w)+i$, for $0\leq i\leq k-1$.

We complete the proof in the following cases.

Case 1: $|R(xwy_1\cdots y_{k-1})|=2$.

When $R(xwy_1\cdots y_{k-1})=\{s,\alpha\}$, by assumptions $y_k=\beta$.
Then $R(xwy)=R(wy)=\{t,r\}$.

When $R(xwy_1\cdots y_{k-1})=\{t,r\}$, by assumptions $y_k=s$.
By Lemma 4.1, Corollary 4.2, $r, t\notin R(xwy)$, hence $R(xwy)=R(wy)=\{s\}$.
It is easy to see that $xwy=(x)(w)(y)$.

Case 2: $R(xwy_1\cdots y_{k-1})=\{y_{k-1}\}$.

We have $R(xwy)\supseteq R(wy)$.
If $R(xwy)=R(wy)$, it is done.

Assume that $R(xwy)\supsetneqq R(wy)$, then $R(xwy)=\{t,r\}$, or $\{s, \alpha\}$.
If $R(xwy)=\{t,r\}$, by the assumption, we get $\{y_{k-1}, y_k\}=\{t, r\}$, hence $R(xwy)=R(wy)=\{t, r\}$, which contradicts the assumption.

If $R(xwy)=\{s,\alpha\}$, by the assumption, we get $\{y_{k-1}, y_k\}=\{s, \alpha\}$.
By 1.1.(d), there exists $u_1\in W$, such that $xwy=(u_1)(w_{s\alpha})$,
write $wy=wy_1\cdots y_i s^a(\alpha s)^b \alpha^c$, here $i$ minimal, such that
$R(y_1\cdots y_i)=\{y_i\}=\{\beta\}$, $a,c=0$ or 1, $a+2b+c<m_{s\alpha}$. Let $u_{s\alpha}=w_{s\alpha}\alpha^c(\alpha s)^{-b}s^a$, $l(u_{s\alpha})\geq 1$. Then $xwy_1\cdots y_i=(u_1)(u_{s\alpha})$.
Hence $i>0$ and $R(wy_1\cdots y_i)=\{y_i\}$ is impossible.
Then $i=0$, or $R(wy_1\cdots y_i)\supsetneqq \{y_i\}$.

If $i=0$, then $y\in W_{s\alpha}$, which contradicts to $s,\alpha$ are not in $L(y)$.

 If $R(wy_1\cdots y_i)\supsetneqq \{y_i\}$ and $i\geq 1$, we have $R(wy_1\cdots y_i)=\{t,r\}$, or $\{s,\beta\}$.
 $xwy_1\cdots y_i=(u_1)(u_{s\alpha})$, when $R(wy_1\cdots y_i)=\{t,r\}$, then $L(u_{s\alpha})=\{\alpha\}$,
 furthermore if $l(u_{s\alpha})=1$, it easy to see that $R(xwy)=R(wy)$, which contradicts the assumption.
 Hence $L(u_{s\alpha})=\{\alpha\}$,
and $l(u_{s\alpha})\geq 2$, by Lemma 4.3, $l(u_{s\alpha})=2$, $u_{s\alpha}=s\alpha$ , if $i\geq 2$, $y_{i-1}=s$, then $s, \alpha \in R(wy_1\cdots y_{i-1})$, then we get equality, $(xwy_1\cdots y_{i-2})(s\beta)=(u_1)(s\alpha)$, which contradicts Lemma 4.4.
If $i=1$, $y_1=\beta$, $xw\beta=(u_1)(s\alpha)$, then $\alpha\in R(xw)=R(w)$, since $l(w)\geq 4$, $(\tilde{x})(s\alpha)(s\beta)=(u_1)(s)$, here $(\tilde{x})(s\alpha s)=xw\alpha$. Then by 1.1.(d), $(\tilde{x})(s\alpha)=(u_2)(s\beta)$, which contradicts Lemma 4.4.

When $R(wy_1\cdots y_i)=\{s,\beta\}$, $L(u_{s\alpha})=\{s\}$, when $l(u_{s\alpha})=1$, it contradicts the assumption.
When $l(u_{s\alpha})\geq 2$, since $xwy_1\cdots y_i=(u_1)(u_{s\alpha})=(x)(u_2)(w_{s\beta})$, which contradicts Lemma 4.3.

Hence $R(xwy)=R(wy)$, and $xwy=(x)(w)(y)$.

\medskip
\textbf{Remark }.
From the prove of Lemma 4.5 we see that if $xw=(x)(w)$, $wy=(w)(y)$, write $y=y_1\cdots y_k$, any reduced decomposition, and $R(xwy_1)\neq \{t,r\}$.  Furthermore if $R(xwy_1\cdots y_i)=R(wy_1\cdots y_i)$, for $i\leq 2$.
Then $i\geq 3$, $R(xwy_1\cdots y_i)=R(wy_1\cdots y_i)$. Hence $xwy=(x)(w)(y)$.

\medskip
\textbf{The Notations,}
from now on, we assume that $m_{sr}\geq 5$ and $m_{st}\geq 4$.

\medskip

\noindent{\bf  Lemma 4.6.}
$x,y\in W$, assume that $t$, $r\notin R(x)\cup L(y)$, then deg $f_{xtr, y , z}\leq 1$, for all $z\in W$.

Proof. We discuss it in two cases.

Case 1: When there is no $x'\in W$, such that $x=(x')(w_{s\alpha}\alpha)$.
Claim that $xtry=(x)(tr)(y)$.
By the Remark above, we only have to check whether $R(xtry)=R(try)$, when $l(y) \leq 2$.

$R(xtr)=R(tr)=\{t,r\}$, it is clear.

$R(xtrs)=R(trs)=\{s\}$, by Lemma 4.1 and Corollary 4.2.

$y=s\alpha$, $R(trs\alpha)=\{\alpha\}$, by Lemma 4.3, $\beta\notin R(xtrs\alpha)$.
Assume that $s\in R(xtrs\alpha)$, then $s\in R(x\beta)$, then there exists $x'\in W$, such that $x=(x')(w_{s\beta}\beta)$, which contradicts the assumption.
Hence $R(xtrs\alpha)=\{\alpha\}$. Hence the claim.

Case 2: When there exists $x'\in W$, such that $x=(x')(w_{s\alpha}\alpha)$, then $xtr=x'w_{s\alpha}\beta$.
If $l(xtry)=l(x)+2+l(y)$, nothing needs to prove. Write $y=y_1\cdots y_k$, reduced decomposition.
Assume that there exists $i<k$, $i$ minimal,
such that $l(w_{s\alpha}\beta y_1\cdots y_i)\leq l(w_{s\alpha}\beta y_1\cdots y_{i-1})$.
By strong exchange condition, we get $s\beta y_1\cdots y_{i-1}=\beta y_1\cdots y_i$,
and $s,\beta\in L(ty_1\cdots y_i)$, write $\beta y=(w_{s\beta})(y')$, by Lemma 1.3.
Since $l(\beta y)=l(y)+1$, $s\beta y_1\cdots y_{i-1}=\beta y_1\cdots y_i$, $s\beta y_1\cdots y_{i-1}y_{i+1}\cdots y_k$ is a reduced expression,
so is $\beta y_1\cdots y_{i-1}y_{i+1}\cdots y_k$. Write $\beta y_1\cdots y_{i-1}y_{i+1}\cdots y_k=(sw_{s\beta})(y')$.
\begin{eqnarray*}
\tt_{w_{s\alpha}\beta}\tt_y & = & \tt_{w_{s\alpha}s}\tt_{s\beta y_1\cdots y_{i-1}}\tt_{y_i}\tt_{y_{i+1}\cdots y_k}\\
               & = & \xi\tt_{w_{s\alpha}s}\tt_{\beta y}+\tt_{w_{s\alpha}s}\tt_{\beta y_1\cdots y_{i-1}y_{i+1}\cdots y_k}\\
               & = & \xi\tt_{w_{s\alpha}s}\tt_{w_{s\beta}y'}+\tt_{w_{s\alpha}s}\tt_{sw_{s\beta}y'}
\end{eqnarray*}

Since $l(w_{s\beta})\geq 4$, $R((x')(w_{s\alpha}s))=L(y')=\{\alpha\}$,
$\xi\tt_{x'w_{s\alpha}s}\tt_{w_{s\beta}y'}=$ $\xi\tt_{x'w_{s\alpha}sw_{s\beta}y'}$, by Lemma 4.5. 

Since $m_{sr}\geq 5$,
$\tt_{x'w_{s\alpha}s}\tt_{sw_{s\beta}y'}=\tt_{(x'w_{s\alpha}s)(sw_{s\beta}y')}$,  by Lemma 4.5.

Hence the lemma is proved.

\medskip

\noindent{\bf  Lemma 4.7.}
$x,y\in W$, when $R(x)=\{\alpha\}$, $L(y)=\{s\}$, then deg $f_{xw_{s\beta}, try , z}\leq 2 $, for all $z\in W$.

Proof.
$$\tt_{xw_{s\beta}}\tt_{try}=\xi\tt_{xw_{s\beta}}\tt_{\alpha y}+\tt_{xw_{s\beta}\beta}\tt_{\alpha y}$$
By Lemma 4.6, we see that deg $\xi\tt_{xw_{s\beta}}\tt_{\alpha y}\leq 2$, since $R((x)(w_{s\beta}\beta))=L(y)=\{s\}$.
As for the part $\tt_{xw_{s\beta}\beta}\tt_{\alpha y}$, when $L(\alpha y)=\{\alpha\}$,
then $w_{s\beta}\beta\alpha y=(w_{s\beta}\beta)(\alpha y)$,
then it is easy to check that $R(xw_{s\beta}\beta)=R(w_{s\beta}\beta)=\{s\}$.
$R(xw_{s\beta}\beta\alpha)=R(w_{s\beta}\beta\alpha)=\{\alpha\}$, it is easy to see that $\beta\notin R(xw_{s\beta}\beta\alpha)$.
If $s\in R(xw_{s\beta}\beta\alpha)$ , then it contradicts Lemma 4.4.
$R(xw_{s\beta}\beta\alpha s)=R(w_{s\beta}\beta\alpha s)$. By Lemma 4.1 or Corollary 4.2, $\beta\notin R(xw_{s\beta}\beta\alpha s)$.
If $\alpha\in R(xw_{s\beta}\beta\alpha s)$ , then it contradicts Lemma 4.3 or Lemma 4.4, since $m_{sr}\geq 5$.
Then by the remark after Lemma 4.6, we have $\tt_{xw_{s\beta}\beta}\tt_{\alpha y}=\tt_{xw_{s\beta}\beta\alpha y}$.

Since $\beta\notin L(\alpha y)$, we have to consider the only left case $L(\alpha y)=\{s,\alpha\}$, write $\alpha y=(w_{s\alpha})(y')$.
\begin{eqnarray*}
\tt_{xw_{s\beta}\beta}\tt_{\alpha y} & = & \tt_{xw_{s\beta}\beta}\tt_{w_{s\alpha}y'}\\
                                    & = &\xi \tt_{(x)(w_{s\beta}\beta s)}\tt_{(w_{s\alpha})(y')}+\tt_{(x)(w_{s\beta}\beta s)}\tt_{(sw_{s\alpha})(y')}
\end{eqnarray*}

Since $s\notin R((x)(w_{s\beta}\beta s))$, there two possibilities.

When $R((x)(w_{s\beta}\beta s))=\{\beta\}$, then by Lemma 4.6, deg $\tt_{xw_{s\beta}\beta}\tt_{\alpha y}\leq 2$.

When $R((x)(w_{s\beta}\beta s))=\{t,r\}$, then $R((x)(w_{s\beta}\beta s\beta))=\{s, \alpha\}$, then $\beta=t$, $m_{st}=4$, since $m_{sr}\geq 5$, write $(x)(s)=(x')(w_{sr})$.

\begin{eqnarray*}
\tt_{xw_{s\beta}\beta}\tt_{\alpha y} & = & \tt_{xw_{s\beta}\beta}\tt_{w_{s\alpha}y'}\\
                                    & = &\xi \tt_{(x)(w_{s\beta}\beta s)}\tt_{(w_{s\alpha})(y')}+\tt_{(x)(w_{s\beta}\beta s)}\tt_{(sw_{s\alpha})(y')}\\
                                    & = & \xi \tt_{(x')(w_{sr}t)}\tt_{(w_{sr})(y')}+\tt_{(x')(w_{sr}t)}\tt_{(sw_{sr})(y')}\\
& = & \xi ^2\tt_{(x')(w_{sr}rt)}\tt_{(w_{sr})(y')}+\xi\tt_{(x')(w_{sr}rt)}\tt_{(rw_{sr})(y')}\\
&    &    +\xi\tt_{(x')(w_{sr}rt)}\tt_{(sw_{sr})(y')}+\tt_{(x')(w_{sr}rt)}\tt_{(rsw_{sr})(y')}
\end{eqnarray*}

Since $m_{sr}\geq 5$, $l(sw_{sr})\geq 4$, $l(rw_{sr})\geq 4$, $R((x')(w_{sr}rt))=L(y')=\{t\}$, by Lemma 4.5,
$$\xi ^2\tt_{(x')(w_{sr}rt)}\tt_{(w_{sr})(y')}=\xi ^2\tt_{(x')(w_{sr}rt)(w_{sr})(y')}$$ $$\xi\tt_{(x')(w_{sr}rt)}\tt_{(rw_{sr})(y')}=\xi\tt_{(x')(w_{sr}rt)(rw_{sr})(y')}$$
$$\xi\tt_{(x')(w_{sr}rt)}\tt_{(sw_{sr})(y')}=\xi\tt_{(x')(w_{sr}rt)(sw_{sr})(y')}$$

Since $l(w_{sr}r)\geq 4$, $R(x')=L(trsw_{sr})=\{t\}$, by Lemma 4.5,
$$\tt_{(x')(w_{sr}rt)}\tt_{(rsw_{sr})(y')}=\tt_{(x')(w_{sr}r)}\tt_{(trsw_{sr})(y')}=\tt_{(x')(w_{sr}rt)(rsw_{sr})(y')}$$

Hence the lemma is proved.

\medskip

\noindent{\bf  Lemma 4.8.}
$x,y\in W$, then deg $f_{xtr, y , z}\leq a$, for all $z\in W$, here $a=$ max $\{m_{sr}, m_{st}\}$.

Proof. We proof the lemma 3 cases.

Case 1:  $R(x)=\{\alpha\}$ and $L(y)=\{\beta\}$.
Write $x=(x')(\alpha)$, $y=(\beta)(y')$, here $R(x')=L(y')=\{s\}$.
Then $\tt_x\tt_y=\tt_{x'tr}\tt_{y'}$, by Lemma 4.6, it is done.

Case 2: $R(x)\bigcup L(y)\neq \{t,r\}$, and furthermore $R(x)\neq \{t, r\}$ or $L(y)\neq \{t,r \}$ .
Let $I=\{s, \alpha\}$, $W_I$ is the parabolic subgroup generated by $I$.
Let $x'$ (resp. $y'$)be the element of minimal length in the coset $xW_I$ (resp.$W_Iy$ ).
Let $w, u\in W_I$ be such that $x=x'w$ and $y=uy'$.
We take proper $\alpha$ here, such that $l(u),l(w)\geq 1$ and $l(w)+l(u)\geq 3$.
Next we use induction on $l(x)+l(y)$, denote $k=l(x)+l(y)$, if $k \leq 2a+1$, nothing needs to prove.

Assume that $k>2a+1$, and the lemma is true for $x'', y''$ with $l(x'')+l(y'')< k$, $R(x'')\neq \{t, r\}$, $L(y'')\neq \{t,r \}$,
$ R(x'')\bigcup L(y'')\neq \{t, r\}$.
$$\tt_x\tt_y=\sum_{v\in W_I}f_{w,u,v}\tt_{x'v}\tt_{y'}$$
When $l(v)\geq 4$, by Lemma 4.5, $\tt_{x'v}\tt_{y'}=\tt_{x'vy'}$.
When $l(v)=0$, since min $\{l(x'), l(y')\}\leq k-1$, by induction hypotheses, we see that the degrees of $f_{x', y', z}$
are not greater than $a$ for any $z\in W$.
Now consider the case $l(v)=1$. If $v=s$, $\xi\tt_{x's}\tt_{y'}=\tt_{x's}\tt_{sy'}-\tt_{x'}\tt_{y'}$, by induction hypotheses,
we see that deg $\tt_{x's}\tt_{sy'}$ and deg $\tt_{x'}\tt_{y'}$ are not greater than $a$.
If $v=\alpha$, write $x'=(x_1)(s\beta)$, $y'=(\beta s)(y_1)$, here $R(x_1s)=L(sy_1)=\{s\}$. $$\xi\tt_{x'v}\tt_{y'}=\xi^2\tt_{x_1str}\tt_{sy_1}+\xi\tt_{x_1s\alpha}\tt_{sy_1}=\xi^2\tt_{x_1str}\tt_{sy_1}+\tt_{x_1s\alpha}\tt_{\alpha sy_1}-\tt_{x_1s}\tt_{sy_1}$$
By Lemma 4.6, we see that deg $\xi^2\tt_{x_1str}\tt_{sy_1}\leq 3$.
By induction hypotheses, we see that deg $(\tt_{x_1s\alpha}\tt_{\alpha sy_1}-\tt_{x_1s}\tt_{sy_1})$ are not greater than $a$.

Hence deg $\tt_{x'v}\tt_{y'}\leq a-1$, when $l(v)=1$.

When $l(v)=2$, $v=s\alpha$ or $\alpha s$.
Only check the case $v=s\alpha$.
$\tt_{x's\alpha}\tt_{y'}=\tt_{x's}\tt_{\alpha y'}$,
by Lemma 4.6 and 4.7, deg $\tt_{x's}\tt_{\alpha y'}\leq 2$.

When $l(v)=3$.
$v=\alpha s \alpha$, then by Lemma 4.6, deg $\tt_{x'\alpha s}\tt_{\alpha y'}\leq 1$.
$v=s\alpha s$, $$\tt_{x's\alpha s}\tt_{y'}=\tt_{x's\alpha s}\tt_{y'}$$
When $R(x's)=L(sy')=\{s\}$,
it is easy to check that $\tt_{x's\alpha s}\tt_{y'}=\tt_{x's\alpha sy'}$.
When $R(x's)$ or $L(sy')$ $=\{s,\beta\}$, then by Lemma 4.6 and 4.7, we have deg $\tt_{x's\alpha s}\tt_{y'}\leq 2$.

Hence the lemma is true in Case 2.

Case 3: When $R(x)=\{t, r\}$ and $L(y)=\{t,r \}$, write $x=(x')(tr)$ and $y=(tr)(y')$.
$$\tt_x\tt_y=\xi^2\tt_{x'tr}\tt_{y'}+\tt_{x't}\tt_{ty'}+\tt_{x'r}tt_{ry'}-\tt_{x'}\tt_{y'}$$
By Lemma 4.6, deg $(\xi^2\tt_{x'tr}\tt_{y'})\leq 3$.
By case 2, deg $\tt_{x't}\tt_{ty'}+\tt_{x'r}tt_{ry'}-\tt_{x'}\tt_{y'}\leq a$.

Hence the lemma is proved.

\medskip

Until now, we see that Theorem 2.1 is proved.

\bibliographystyle{unsrt}

\end{document}